\documentclass[preprint,3p,11pt,sort&compress]{elsarticle}

 \usepackage{hyperref}
 \hypersetup{backref, colorlinks=true, linkcolor=blue}

\usepackage{stfloats}
\usepackage{subcaption}

 \usepackage{pict2e}

\usepackage{graphicx} % allows inclusion of graphics
\usepackage{booktabs} % nice rules (thick lines) for tables
\usepackage{microtype} % improves typography for PDF
\usepackage{amssymb,amsmath}
\usepackage[ruled,vlined]{algorithm2e}

\begin{document}
\begin{center}
{\bf \Large Convergence of Multi-Level Hybrid Monte Carlo Methods \\
	\vspace{0.1cm}
	for 1-D Particle Transport Problems}
\end{center}

\author[ncsu]{Vincent N. Novellino}
\ead{vnnovell@ncsu.edu}
\author[ncsu]{Dmitriy Y. Anistratov}
\ead{anistratov@ncsu.edu}

\address[ncsu]{Department of Nuclear Engineering,
North Carolina State University Raleigh, NC}

\begin{frontmatter}
\begin{abstract}
We present in this paper a hybrid, Multi-Level Monte Carlo (MLMC) method for solving the neutral particle transport equation.
MLMC methods, originally developed to solve parametric integration problems, work by using a cheap, low fidelity solution as a base solution and then solves for additive correction factors on a sequence of computational grids.
The proposed algorithm works by generating a scalar flux sample using a Hybrid Monte Carlo method based on the low-order Quasidiffusion equations.
We generate an initial number of samples on each grid and then calculate the optimal number of samples to perform on each level using MLMC theory.
Computational results are shown for a 1-D slab model to demonstrate the weak convergence of considered functionals.
The analyzed functionals are integrals of the scalar flux solution over either the whole domain or over a specific subregion.
We observe the variance of the correction factors decreases faster than increase in the cost of generating a MLMC sample grows. 
The variance and costs of the MLMC solution are driven by the coarse grid calculations. 
Therefore, we should be able to add additional computational levels at minimal cost since fewer samples would be needed to converge estimates of the correction factors on subsequent levels.
\end{abstract}
\end{frontmatter}

	\section{Introduction}\label{sec:1}
	
	Monte Carlo (MC) methods for solving the particle transport problem involve simulating the movement of particles through a problem domain and collecting tallies of "events" on spatial, energy, and temporal grids.
	The quality of the solution depends on the number of particles simulated that travel to a particular grid element. 
	Estimators for the mean will converge by the Central Limit Theorem $\frac{\big<{X}_N\big> - \mu}{\sigma\sqrt{N}} \xrightarrow{d} \mathcal{N}(0,1)$, where $\big<{X}_N\big> = \frac{1}{n}\sum_{n=1}^N X(\omega_n)$, $X$ is our random variable of interest, and $\omega_n \in \Omega$ is a randomly sampled event from the event space $\Omega$. 
	In applications where the global solution for the particle transport problem is needed, MC methods can be slow to converge, especially so for large problems with small grid elements that are in low probability regions.
	If we increase the size of our grid elements then we can, on average, increase the number of particles that contribute tallies to the estimator in that element, though we decrease the overall resolution of the solution.
	What if there is a way to restore the resolution without adding significant costs to the calculation?
	Let us introduce the idea of multi-level Monte Carlo (MLMC) methods.
	
	MLMC methods were initially developed for solving parametric integration problems \cite{heinrich-LSSC2001} but extensions to solving stochastic differential and partial differential equations have been made \cite{giles-AN-2015}.
	The basic idea of MLMC is to use a cheap, low fidelity solution that captures the essential structure of the solution and then solve for additive correction terms that yields the solution on a higher fidelity target grid.
	Then, using estimators of the mean, variance, and computational costs (run times), the MLMC algorithm optimizes the amount of work to be performed on each grid to yields a lower variance and cheaper solution in comparison to a standard MC calculation on the target grid.
	Early optimization algorithms considered a scalar functional as a quantity of interest, i.e. hydraulic conductivity for groundwater flow simulations \cite{cliffe-2011}, however one can apply MLMC theory to a vector of output functionals by checking convergence of each component \cite{giles-AN-2015}.
	MLMC optimization depends on the coarse grid solution being computationally less expensive than the finer grid solutions and the variance of the correction terms should decrease with increasing grid fidelity.
	The algorithm moves the majority of the computational effort to the coarsest level, we have smaller computational costs relative to an equivalent MC simulation on the target grid.
	In addition, one can generate a low variance solution on a coarse grid cheaper than on the fine grid, meaning the overall uncertainty in the solution is reduced using MLMC.

	Hybrid MC methods (HMC) have been developed for the global transport problem by solving for moments of the angular flux with closures estimated using MC.
	These previous methods primarily sought to accelerate fission source convergence in criticality calculations \cite{larsen-yang-nse-2008,Lee-Joo-Lee-Smith,wolters-nse-2013,Willert-Knoll-Kelley-Park} or to remove effective scattering events in Implicit MC calculations \cite{pozulp-mc2023}.
	This is a group of HMC methods that have been developed for solving the particle transport problem, but for brevity we discuss those that most related to the work we will present in this paper.
	
	In this work, we apply the MLMC algorithm to a HMC method that uses low-order equations to solve for the zeroth moment of the angular flux in 1-D.
	The HMC method we utilize estimates closure terms for the low-order moment equations discretized by a finite volume scheme \cite{vnn-dya-ans-annual-2024}.
	For each computational level, except for the coarsest, we solve the low-order quasidiffusion (LOQD) equations on a coarse and fine grid. 
	Then the solution functional and correction for the functional are evaluated using the HMC results.
	The functional and correction results are used in the MLMC optimization portion of the algorithm to decide how many HMC solves are needed per level.
	
	We organize the remainder of this work by first introducing the relevant low-order moment equation and the finite volume discretization at the beginning of Section \ref{sec:2}.
	Then we will show how we will use Monte Carlo to estimate the LOQD closures in Section \ref{sec:hmc}.
	Section \ref{sec:hmlmc} will describe the hybrid MLMC method in detail, giving the reader a quick overview of the relevant theory, and then showing the optimization algorithm we implemented.
	The results for a selection of test problems with increasing scattering ratios are then shown in Section \ref{sec:res} to demonstrate the convergence of the method.
	Finally, conclusions and future work considerations will be discussed in Section \ref{sec:con}.

	\section{1-D Particle Transport} \label{sec:2}
	
	We consider the 1-D, steady state particle transport equation with isotropic scattering and source:
	\begin{equation*} \label{t-eq}
		\mu \frac{\partial \psi}{\partial x}(x,\mu) + \Sigma_t(x) \psi(x,\mu) = \frac{\Sigma_s(x)}{2} \int_{-1}^1 \psi(x,\mu') d\mu' + \frac{Q(x)}{2} \, ,
	\end{equation*}
	\[x\in D, \quad D = [0,X], \quad \mu \in [-1,1] 
	\]
	\begin{equation*}
		\psi(0,\mu) = 0 , \ \mu > 0 \, , \quad 
		\psi(L,\mu) = 0 , \ \mu < 0 \, .
	\end{equation*}
	The second order form the LOQD equations are given by \cite{gol'din-cmmp-1964}:
	\begin{equation*}
		-\frac{d}{dx} \frac{1}{\Sigma_t(x)} \frac{d(E\phi)}{dx}(x)  + \Sigma_a(x) \phi(x) = Q(x) \, ,
	\end{equation*}
	where the QD (Eddington) factor is defined through
	\begin{equation*} \label{QDf}
		E(x) = \int_{-1}^{1} \mu^2 \psi(x,\mu) d \mu \bigg/  \int_{-1}^{1} \psi(x,\mu) d \mu \, .
	\end{equation*} 
	Boundary factors for the LOQD method can be defined by:
	\begin{equation*}
		B_0 = \frac{ \int_{-1}^{0}\mu \psi(0,\mu) d\mu}
		{\int_{-1}^{0} \psi(0,\mu) d\mu}\; , \quad
		B_X= \frac{ \int_{0}^{1}
			\mu \psi(X,\mu) d\mu}
		{\int_{0}^{1} \psi(X,\mu) d\mu}\; .
		\label{eqn:qd_bf}
	\end{equation*}

	\subsection{Hybrid MC} \label{sec:hmc}
	
	To derive our Hybrid MC method let us first consider the finite volume discretized version of the LOQD equations \cite{dya-vyag-ttsp}:
	\begin{equation*}
		-  \frac{ E_{i+1}\phi_{i+1} - E_i \phi_i  }{\Sigma_{t,i+1/2} \Delta x_{i+1/2}} 
		+ \frac{ E_i\phi_i - E_{i-1} \phi_{i-1}  }{\Sigma_{t,i-1/2} \Delta x_{i-1/2}}    +
		\Sigma_{a,i} \Delta x_i \phi_i =  Q_i \Delta x_i \, ,
	\end{equation*}
	where $i$ is the spatial cell index $(i=1, \dots, I)$, the integer subscript indicates a quantity averaged over the $i^{th}$ cell $g_i = [x_{i-1/2}, x_{i+1/2}]$,  $\Delta x_i = x_{i+1/2} - x_{i-1/2}$,
	\begin{equation*}
		\Sigma_{t,i+1/2} =\frac{\Sigma_{t,i}\Delta x_i  + \Sigma_{t,i+1}\Delta x_{i+1}}{\Delta x_i  + \Delta x_{i+1}} \, , \quad
		\Delta x_{i+1/2}=\frac{1}{2}(\Delta x_i+\Delta x_{i+1}) \, .
	\end{equation*}
	To calculate our cell-average QD factors we use track-length tallies
	\begin{equation*}\label{eqn:pathlength}
		E_i  \! =  \! \frac{\sum_{n=1}^N \! \! \sum_{m=1}^{M_n} \mu_{n,m}^2  w_{n,m} \ell_{n,m}}
		{\sum_{n=1}^N \! \! \sum_{m=1}^{M_n}   w_{n,m} \ell_{n,m}},
	\end{equation*}
	where we have $N$ source particles, $M_n$ tracks, $\mu_{n,m}$ is the cosine of the direction of flight, $w_{n,m}$ is the particle weight, and $l_{n,m}$ is the track-length.
	The QD boundary factors can be estimated on the problem boundary by:
	\begin{equation*}
		B_{0} = \frac{\sum_{n=1}^N \sum_{m=1}^{M_n} w_{n,m} f(\mu_{n,m})}{\sum_{n=1}^N \sum_{m=1}^{M_n} w_{n,m}      \frac{f(\mu_{n,m})}{\mu_{n,m}}}, \,
		B_{X} = \frac{\sum_{n=1}^N \sum_{m=1}^{M_n} w_{n,m} f(-\mu_{n,m})}{\sum_{n=1}^N \sum_{m=1}^{M_{n,m}} w_{n,m} \frac{f(\mu_{n,m})}{\mu_{n,m}}},
	\end{equation*}
	where $f(\mu)$ is a piece-wise function of the form:
	\begin{equation*}
		f(\mu) = 
		\begin{cases}
			1, \quad \mu <0 \\
			0, \quad \mu \geq 0. \\
		\end{cases}
	\end{equation*}

	\subsection{Hybrid Multi-Level Monte Carlo} \label{sec:hmlmc}
	
	We will give a short summary of important results from MLMC theory here, borrowing the notations from \cite{giles-AN-2015}.
	Let us consider a collection of grids, $G_\ell=\{ g_i^\ell\}_{i=1}^{I_\ell}$ for $\ell = 0, 1, \dots, L$ such that $G_0 \subset G_1 \subset G_2 ... \subset G_L$.
	$G_0$ is our coarsest grid and $G_L$ is the finest grid.
	Let $P$ be a functional of the zeroth moment of the angular flux, $\phi(x)$.
	MLMC works by reducing the absolute error in the functional $P$, and if we define $P$ over a specific region of interest, or over the whole domain, then we can control the quality of the solution in that region.
	Let $P_\ell$ be the approximation of the functional on grid $G_\ell$.
	We consider the following estimator of $P_\ell$
	\begin{equation*}
		\big<P_\ell\big> = N_\ell^{-1} \sum_{n=1}^{N_\ell} P_\ell(\omega_{n,\ell}),
	\end{equation*}
	where $\omega_{n,\ell}$ is a sample coming from the probability space $(\Omega, \mathcal{F},P)$.
	In our context, we consider the results of an HMC simulation to be our sample, i.e. we simulate $K$ particle histories, calculate our low-order functionals, and solve the LOQD equation with this set of cell-average QD factors. 
	This provides one realization of the hybrid low-order problem yielding a sample solution.
	Our approximation of $\big<P_L\big>$ using MLMC is given by:
	\[\big<\tilde{P}_L\big> = \big<P_0\big> + \sum_{\ell=1}^L N_{L}^{-1} \sum_{n=1}^{N_{L}} \biggl\{ P_{\ell}(\omega_{n,\ell}) - P_{\ell -1}(\omega_{n,\ell}) \biggr\}, 
 \]
	where we use the same random sample $\omega_{n,\ell}$ to approximate the functional on the fine $G_\ell$ and the coarse $G_{\ell-1}$ for each level.
	Let's define some notations: $\Delta P_{\ell,n} = \Delta P_{\ell}(\omega_{n,\ell}) = P_{\ell}(\omega_{n,\ell}) - P_{\ell -1}(\omega_{n,\ell})$ and $P_{\ell,n} = P_{\ell}(\omega_{n,\ell})$.
	If the cost of one sample on $G_\ell$ is $C_\ell$ and the variance is $V_\ell$, then:
	\[C = \sum_{\ell=1}^{L} N_\ell C_\ell, \quad V = \sum_{\ell=0}^L N_\ell^{-1}V_\ell.  \]
	Given a target $\epsilon^2$ for the variance, we want to minimize the computational costs of the calculations.
	Consider the Lagrangian function:
	\[\mathcal{L}(N_\ell,\lambda) = C(N_\ell) + \lambda^2 V(N_\ell) = \sum_{\ell=0}^{L} N_\ell C_\ell + \sum_{\ell=0}^L \lambda^2 N_{\ell}^{-1}V_\ell.  \]
	Solving for the optimal $N_\ell$ yields $N_\ell = \lambda \sqrt{V_\ell / C_\ell}$.
	Using our target variance of $\epsilon^2$ we can solve for $\lambda = \epsilon^{-2} \sum_{\ell=0}^L \sqrt{V_\ell C_\ell}$.
	Then the total cost is then:
	\[C = \epsilon^{-2} \biggl(\sum_{\ell=0}^L \sqrt{V_l C_l}\biggr)^2. \]
	We can consider three different scenarios: $V_\ell C_\ell$ increasing with $\ell$ meaning the cost is increasing faster than the variance decreases, $V_\ell C_\ell$ decreasing with $\ell$ meaning our variance is increasing faster than the costs, or $V_\ell C_\ell$ staying the same.
	In our analysis we will show we are in the second situation, meaning that our dominant contribution to the cost comes from the first level and $C \approx \epsilon^{-2} C_0 V_0$. 
	We also will consider in this work $G_\ell$ as a geometric sequence where our common ratio is $2$.
	Let's define $\Delta P_\ell$ as the multi-level estimator of $P$ such that
	\[\big<\Delta P_\ell\big> = N_\ell^{-1} \sum_{n=1}^{N_\ell} \bigl\{P_{\ell}(\omega_{n,\ell}) - P_{\ell-1}(\omega_{n,\ell})\bigr\}, \]
	where $P_{-1} = 0$.
	In addition, if $P_{\ell}$ and $P_{\ell-1}$ both approximate $P$, then the difference should decrease as $\ell \to \infty$.
	$\mathbb{V}(\Delta P_\ell) \to 0$ as $\ell \to \infty$.
	If the variance decreases as we add more levels then fewer samples are needed to be simulated for the subsequent levels.
	
	The following theorem summarizes these requirements, plus a few others, and serves as the basis of MLMC theory \cite{giles-AN-2015}:
	Let P be a random variable and $P_\ell$ the approximation on the grid $G_\ell$.
	If $\big<\Delta P_\ell\big>$ is an independent estimator generated by $N_\ell$ samples with expected cost of $C_\ell$ and variance $V_\ell$, positive constants $\alpha$, $\beta$, $\gamma$, $c_i$ for $i\in\{1,2,3\}$, $\alpha \geq \frac{1}{2}\min(\beta,\gamma)$, and conditions 1-4 are met,
	\begin{enumerate}
		\item $|\mathbb{E}(P_\ell - P)| \leq c_1 2^{-\alpha \ell}$
		\item $\mathbb{E}\Big(\big<\Delta P_\ell\big>\Big) = \begin{cases}
			\mathbb{E}(P_0) & \ell = 0 \\
			\mathbb{E}(P_\ell - P_{\ell -1}) & \ell > 0 
		\end{cases}$
		\item $V_\ell \leq c_2 2^{-\beta \ell}$ 
		\item $C_\ell \leq c_3 2^{\gamma \ell}$ 
		\end{enumerate}
	Then there exist a constant $c_4$ so that for any $\epsilon < 1/e$ we can find an $L$ and $N_\ell$ so that:
	\[ \big<\tilde{P}_L\big> = \sum_{\ell=0}^L \big<\Delta P_\ell\big>, \quad \mathbb{E}\Big( \big<\tilde{P}_L\big> - \mathbb{E}(P) \Big)^2 < \epsilon^2, \]
	and the computational complexity is bounded by:
	\[ \mathbb{E}(C) \leq \begin{cases}
		c_4 \epsilon^{-2} & \beta > \gamma \\
		c_4 \epsilon^{-2} \big(\log{\epsilon}\big)^2 & \beta = \gamma \\
		c_4 \epsilon^{-2 - (\gamma-\beta)/\alpha} & \beta < \gamma.
	\end{cases}
	\]

	\subsubsection{MLMC algorithm}
	
	The MLMC scheme is presented in Algorithm \ref{alg:hybrid_mlmc}.
	We consider $L$ levels and check for weak convergence as a post processing step.
	The MLMC algorithm uses $K_{\ell}$ particle histories for generating the grid function of Eddington factors in each spatial cell and boundary factors to form a realization of hybrid discretized LOQD equation the numerical solution of which yields one sample vector of hybrid scalar flux. 
	At the $\ell$-level, we perform  $N_{\ell}$ realization to compute $\big<P_{\ell}\big>$.
	\begin{algorithm}
		\DontPrintSemicolon
		\caption{MLMC Algorithm for 1-D Particle Transport}\label{alg:hybrid_mlmc}
	\small
		$N_\ell \gets 10$ \;
		\ForEach{$u=1, \dots, 3$}{
			
			\ForEach{$\ell = 0, \dots, L$}
			{
				\If{$u> 1$}
				{
					$N_\ell\gets \max\Big(\text{Integer}\big(2\epsilon^{-2} \sqrt{V_\ell / C_\ell} \sum_{\ell=0}^{L} \sqrt{V_\ell C_\ell}\big) - N_\ell, 0\Big)$
					\vspace{-0.1cm}
				}
				\eIf{$\ell = 0$}
				{        
					\ForEach{$n = 1,\dots,N_\ell$}
					{
						Peform MC simulation to estimate Eddington and closure functionals \;
						Perform LOQD solve on grid $G^0$\;
						Evaluate and store $P_\ell$\;
											\vspace{-0.1cm}
					}
					Evaluate $V_\ell = \mathbb{V}(P_\ell)$ and $C_\ell = N_\ell^{-1} \sum_{n=1}^{N_\ell} C_{\ell,n}$ \;
				}
				{
					\ForEach{$n = 1,\dots,N_\ell$}
					{
						Peform MC simulation to estimate Eddington and closure functionals \;                
						Evaluate $V_\ell = \mathbb{V}(P_\ell)$ and $C_\ell = N_\ell^{-1} \sum_{n=1}^{N_\ell} C_{\ell,n}$ \;
						
						Perform LOQD solve on grids $G^{\ell}$ and $G^{\ell-1}$\;
						Evaluate and store $\big<P_\ell\big>$ and $\big<\Delta P_{\ell} \big>$ \;
					}
					Evaluate $V_\ell = \mathbb{V}(\Delta P_{\ell})$ and $C_\ell = N_\ell^{-1} \sum_{n=1}^{N_\ell} C_{\ell,n}$ \;
				}
			}
		}
	\end{algorithm}
	In addition, we examine the effect of updating the estimates of $V_\ell$ and $C_\ell$ for estimating the optimal number of samples per level in this work.
	In this study, we consider  the functional $P_\ell$ as a integral of the scalar flux  over a
	spatial interval $A$
	\[P_\ell = \int_{A} \phi(x) dx \]
	where $A$ is either the whole spatial domain ($A=D$)  or  a cell on the coarsest grid $G_0$ and hence  $A=g_i^0$.
	After completing a MLMC solve by Algorithm \ref{alg:hybrid_mlmc}, we check for weak convergence by:
	\[W_{\hat{\ell}} = \frac{\big<\Delta P_{\hat{\ell}}\big>}{2^{\alpha} - 1} < \frac{\epsilon}{\sqrt{2}}, \quad \mbox{where} \quad \hat{\ell} = \{L-2,L-1,L\}.  \] 
	
	\vspace{-0.2cm}
	\section{Results} \label{sec:res}
	
	The model we consider is a 1-D two material slab of length $X = 1$, $Q(x) = 1$, $\Sigma_t(x) = 1.0$, $c_i =  \Sigma_s/\Sigma_t$, $c_1 = 0.9$, $c_2 = 0.1, 0.5, 0.9$.
	Region 1 is from $0 < x <0.5$ and Region 2 is then $0.5 < x < 1.0$.
	We uniformly mesh the domain for each computational level $G_\ell$ with $I_\ell$ cells, where $I_\ell = 2 I_{\ell -1}$, $I_0 = 16$, and $L =3$.
	Multiple $\epsilon$ values were considered to analyze how the MLMC algorithm converges and distributes work among the computational levels under different convergence criteria.
	We also consider MC samples consisting $K_{\ell}=10^3, 10^4$ ($\ell=0,\ldots,L$) particle histories and we use implicit capture to reduce the variance of the Eddington and boundary factor estimators.
	Changing the number of particles making up a sample should influence the algorithm since it changes the variance of the Eddington and boundary factors, thereby changing the amount of noise in the hybrid solution.
	
\begin{figure*}[t!]
		\centering
		\begin{subfigure}[b]{0.495\textwidth}
			\centering
			\includegraphics[width=\textwidth]{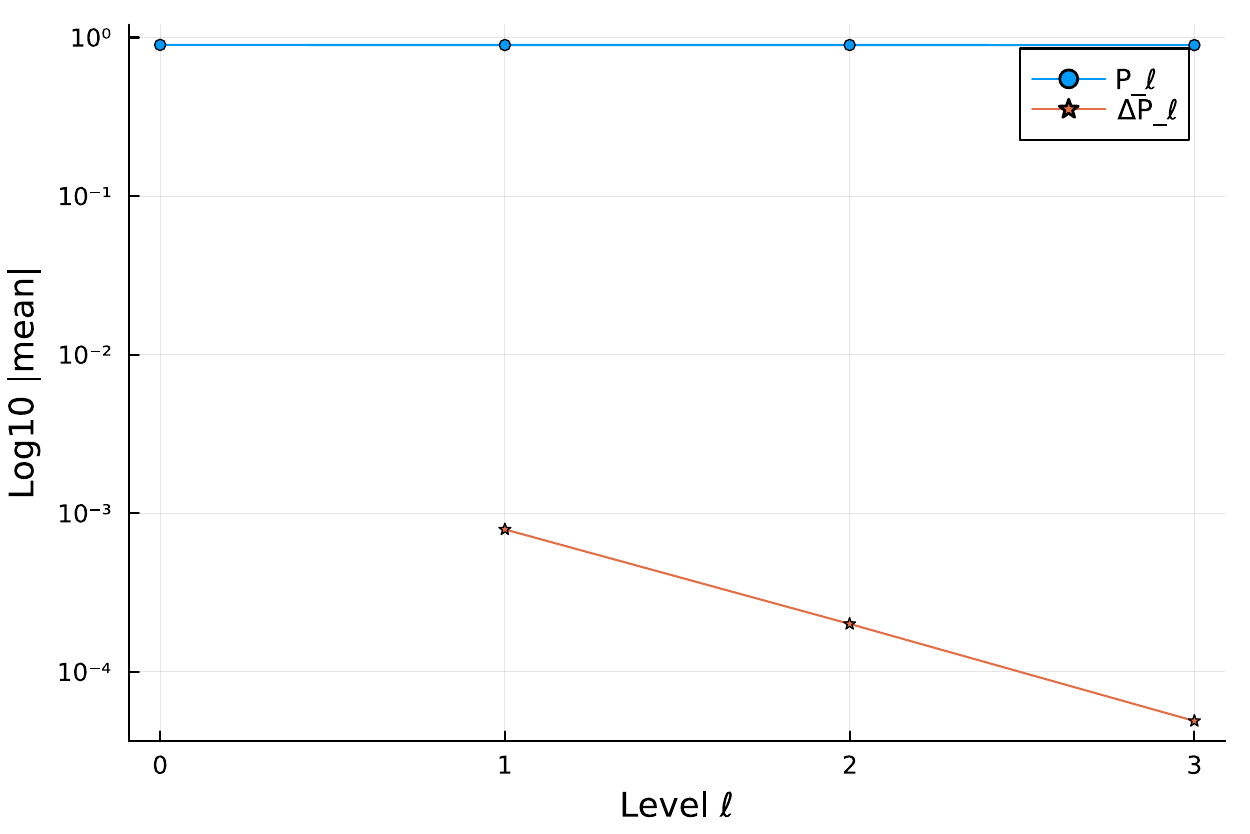}
			\caption*{MLMC Mean Value of Functional}
		\end{subfigure}
		\begin{subfigure}[b]{0.495\textwidth}
			\centering
			\includegraphics[width=\textwidth]{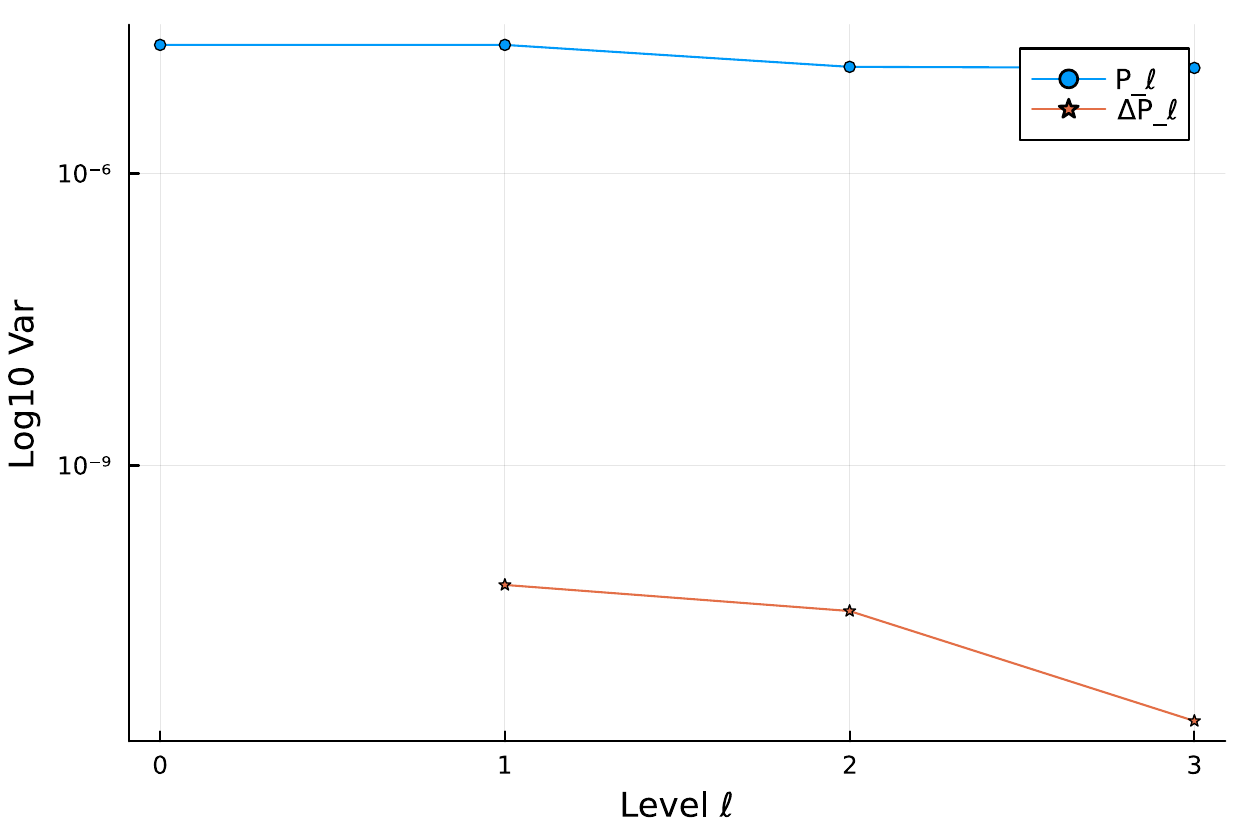}
			\caption*{MLMC Variance of Functional}
		\end{subfigure}
		\begin{subfigure}[b]{0.495\textwidth}
			\centering
			\includegraphics[width=\textwidth]{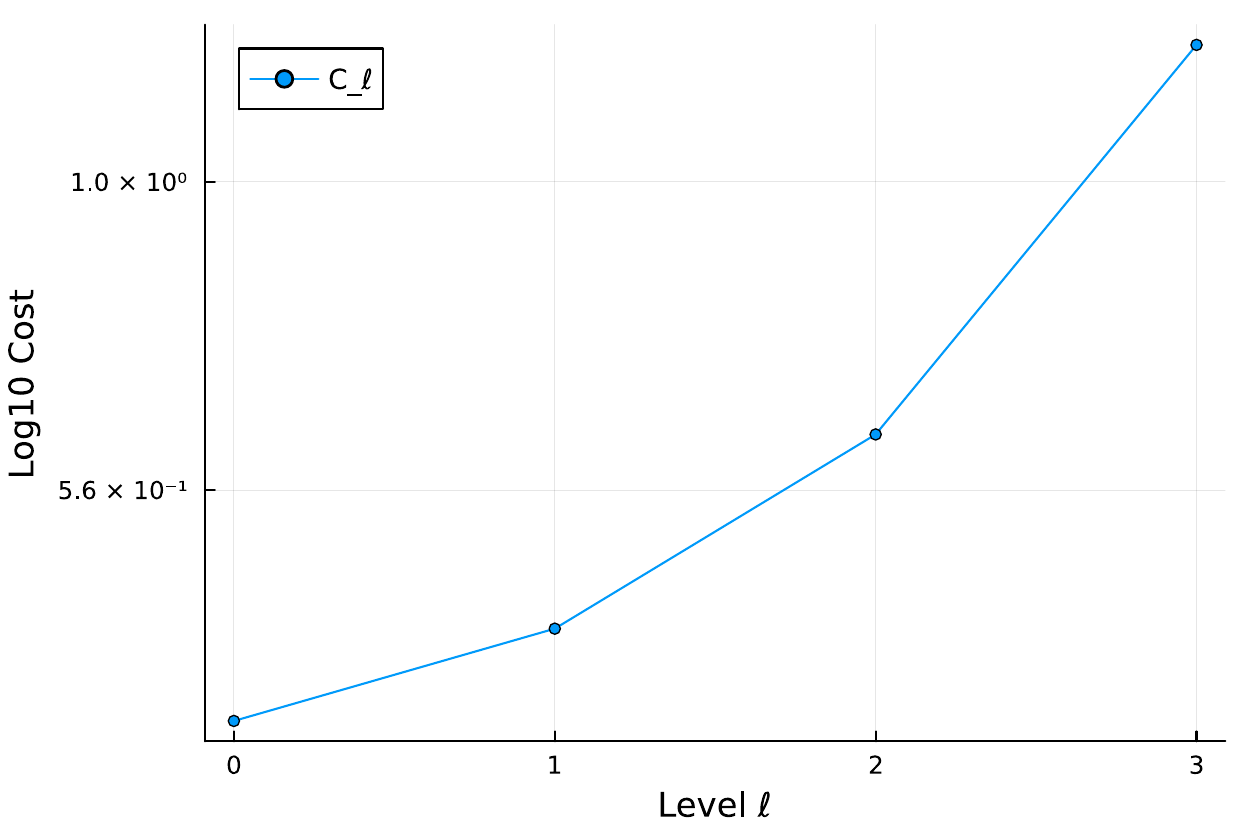}
			\caption*{MLMC Cost Estimates}
		\end{subfigure}
		\begin{subfigure}[b]{0.495\textwidth}
			\centering
			\includegraphics[width=\textwidth]{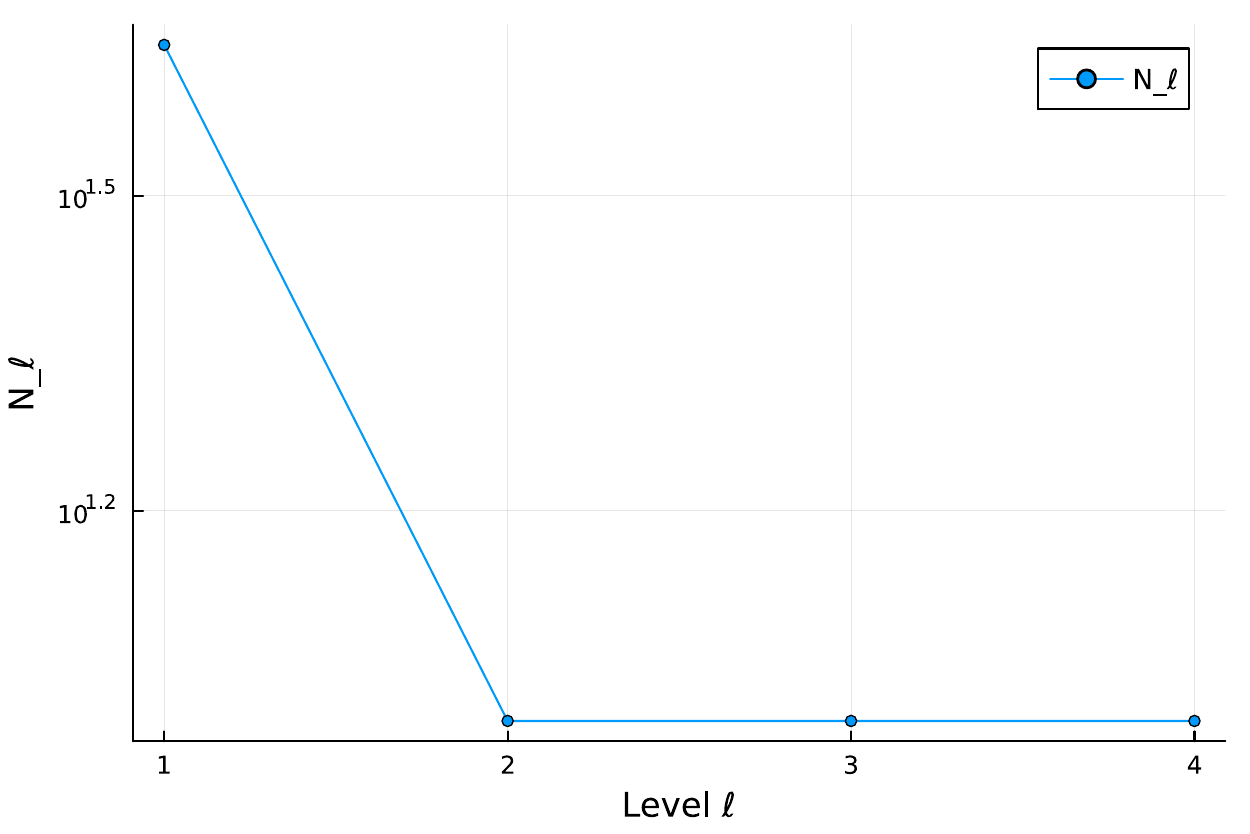}
			\caption*{MLMC Number of Realizations}
		\end{subfigure}
		\begin{subfigure}[b]{0.495\textwidth}
			\centering
			\includegraphics[width=\textwidth]{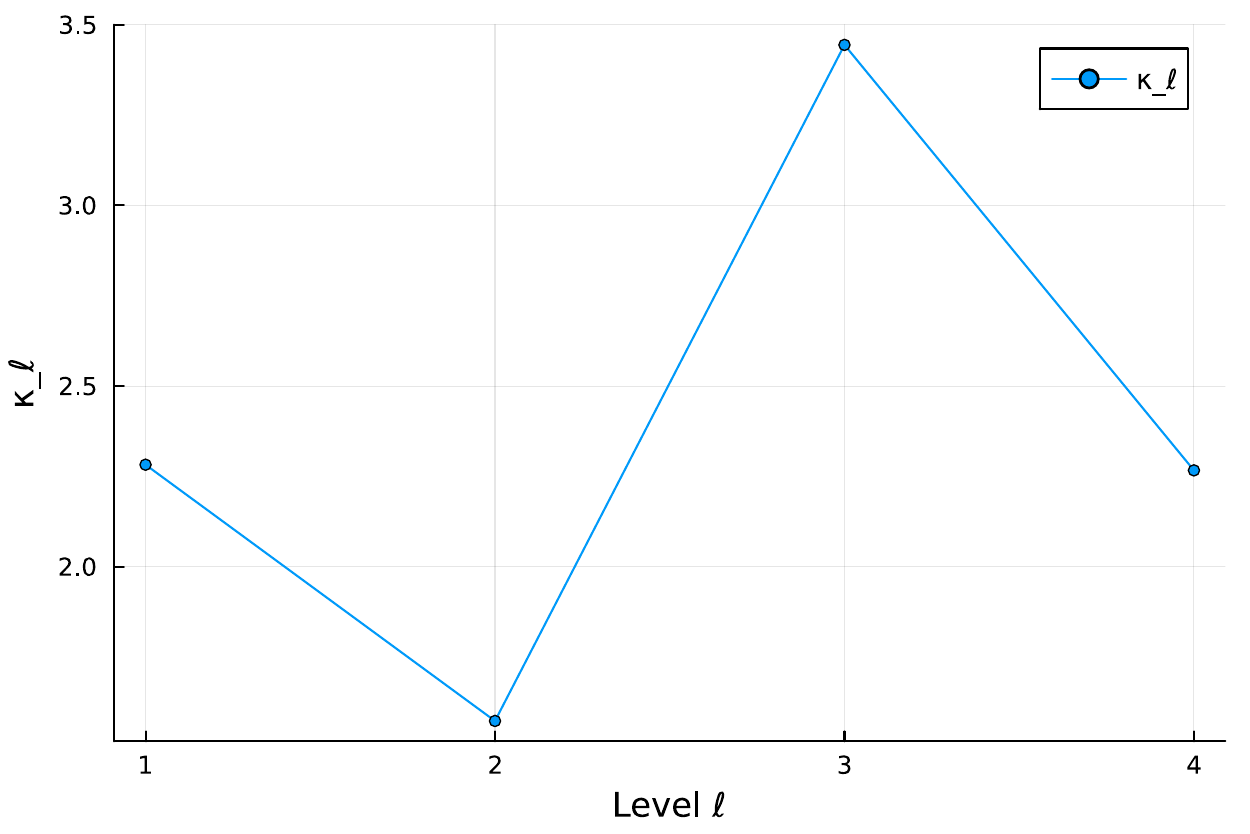}
			\caption*{MLMC Kurtosis}
		\end{subfigure}
		\begin{subfigure}[b]{0.495\textwidth}
			\centering
			\includegraphics[width=\textwidth]{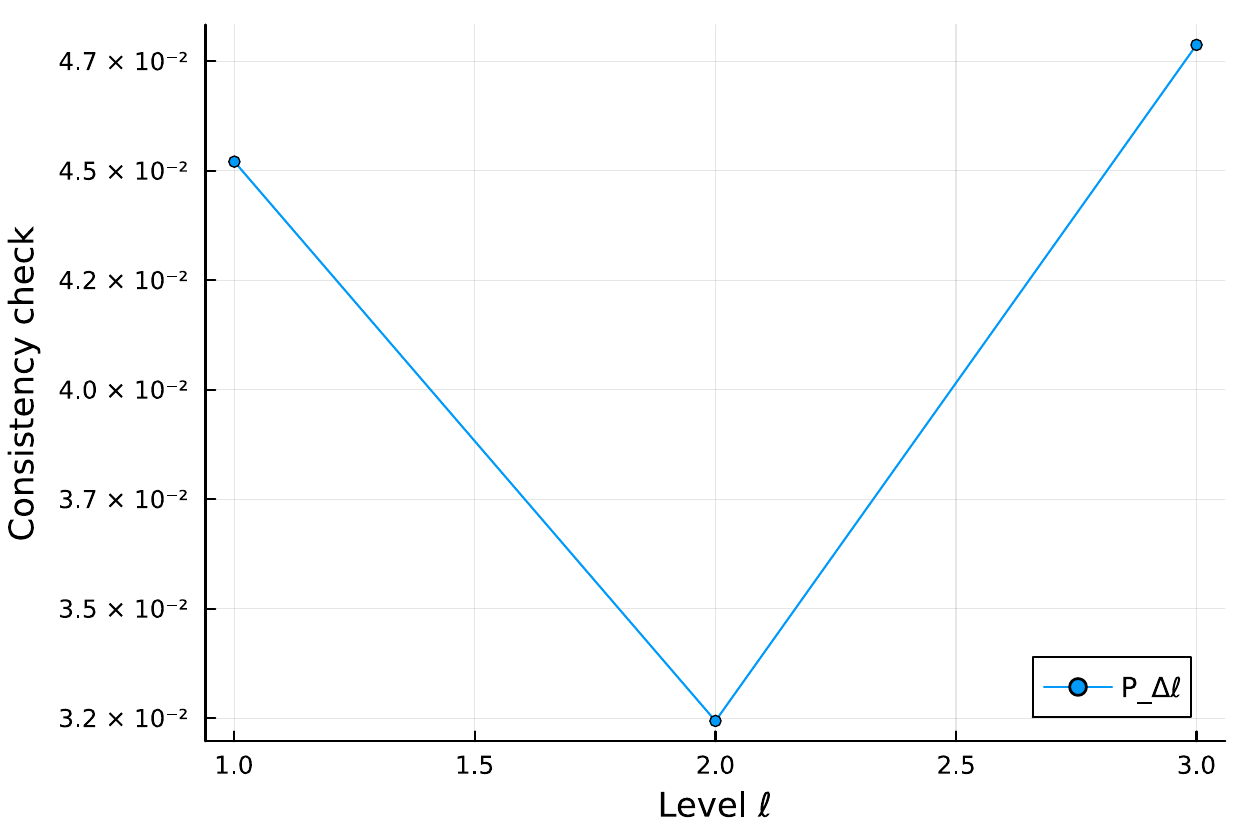}
			\caption*{MLMC consistency check}
		\end{subfigure}
		\caption{MLMC functional results for $c_2 = 0.5$, $K = 10^4$ particle histories, and $\epsilon = 1 \cdot 10^{-3}$}
		\label{fig:mlmc_results_0.5}
	\end{figure*}

\begin{figure*}[t!]
		\centering
		\begin{subfigure}[b]{0.495\textwidth}
			\centering
			\includegraphics[width=\textwidth]{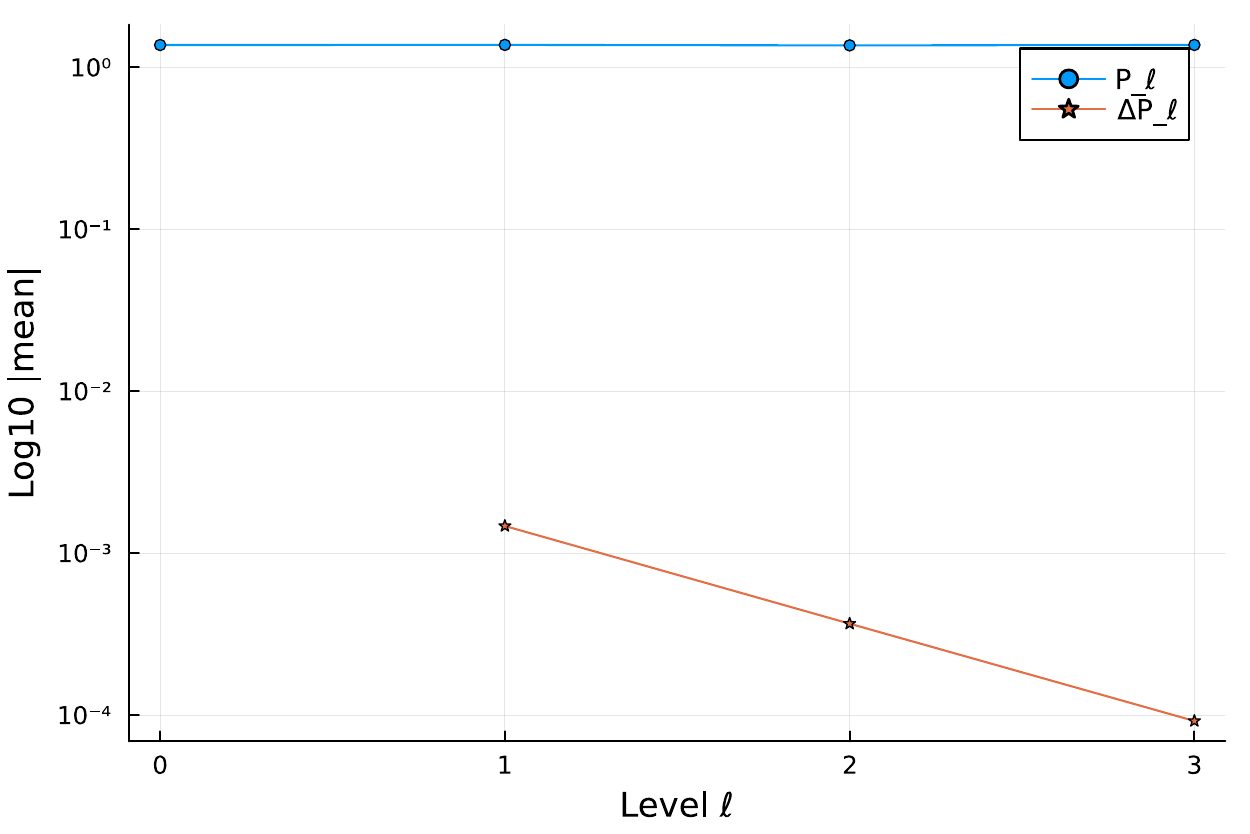}
			\caption*{MLMC Mean Value of Functional}
		\end{subfigure}
		\begin{subfigure}[b]{0.495\textwidth}
			\centering
			\includegraphics[width=\textwidth]{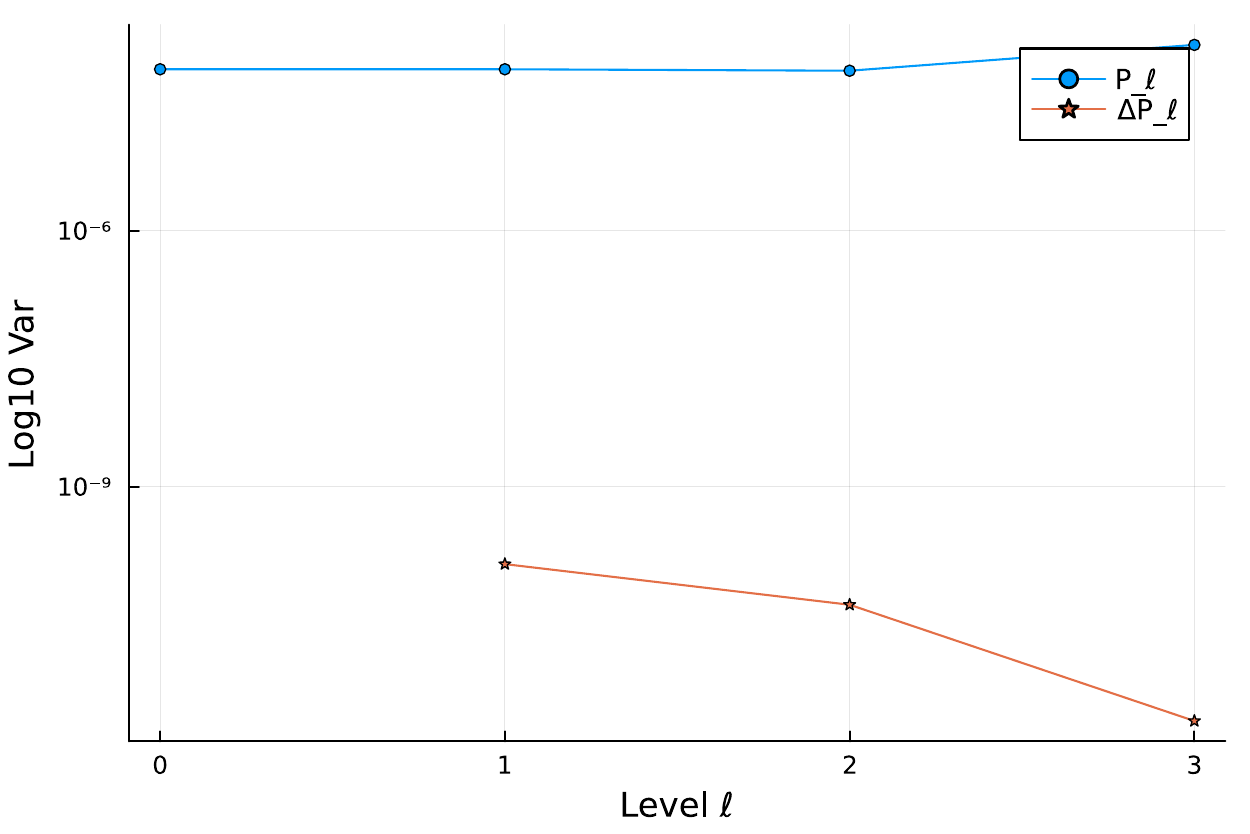}
			\caption*{MLMC Variance of Functional}
		\end{subfigure}
		\begin{subfigure}[b]{0.495\textwidth}
			\centering
			\includegraphics[width=\textwidth]{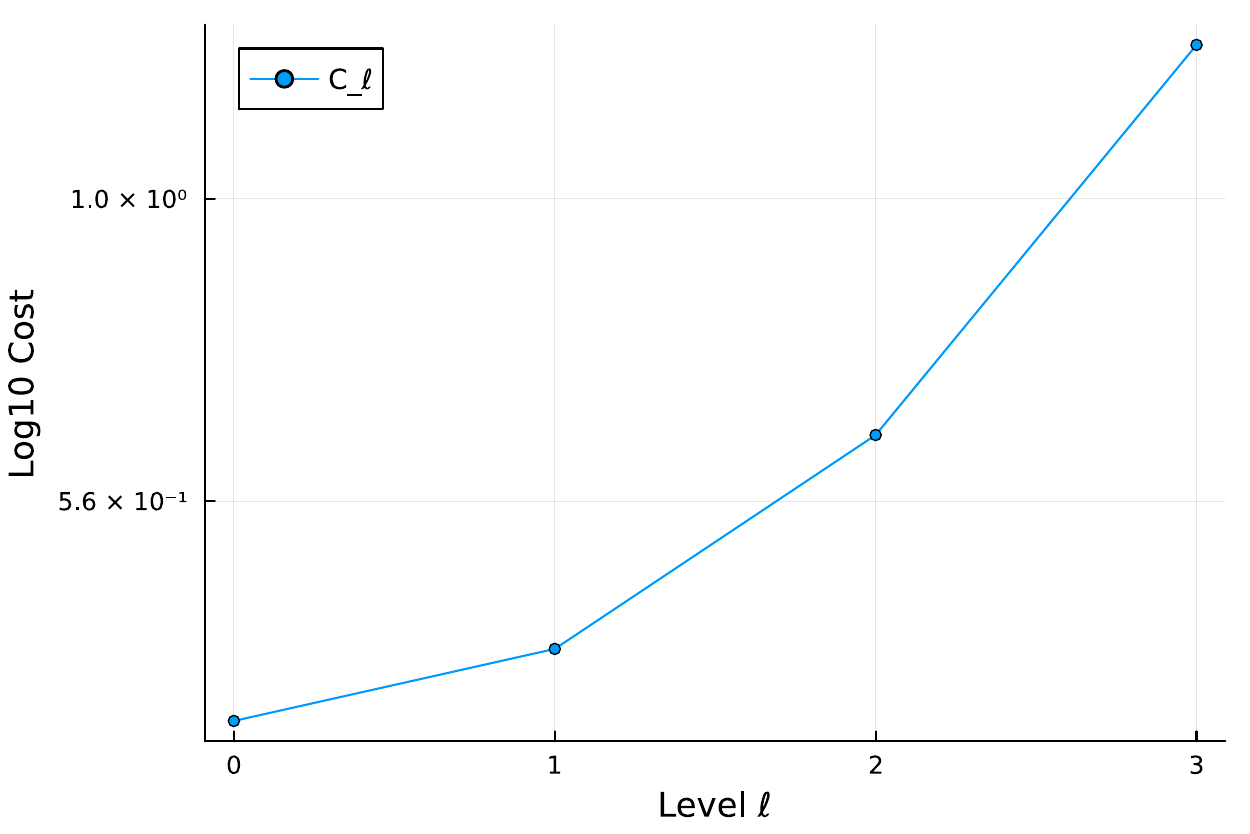}
			\caption*{MLMC Cost Estimates}
		\end{subfigure}
		\begin{subfigure}[b]{0.495\textwidth}
			\centering
			\includegraphics[width=\textwidth]{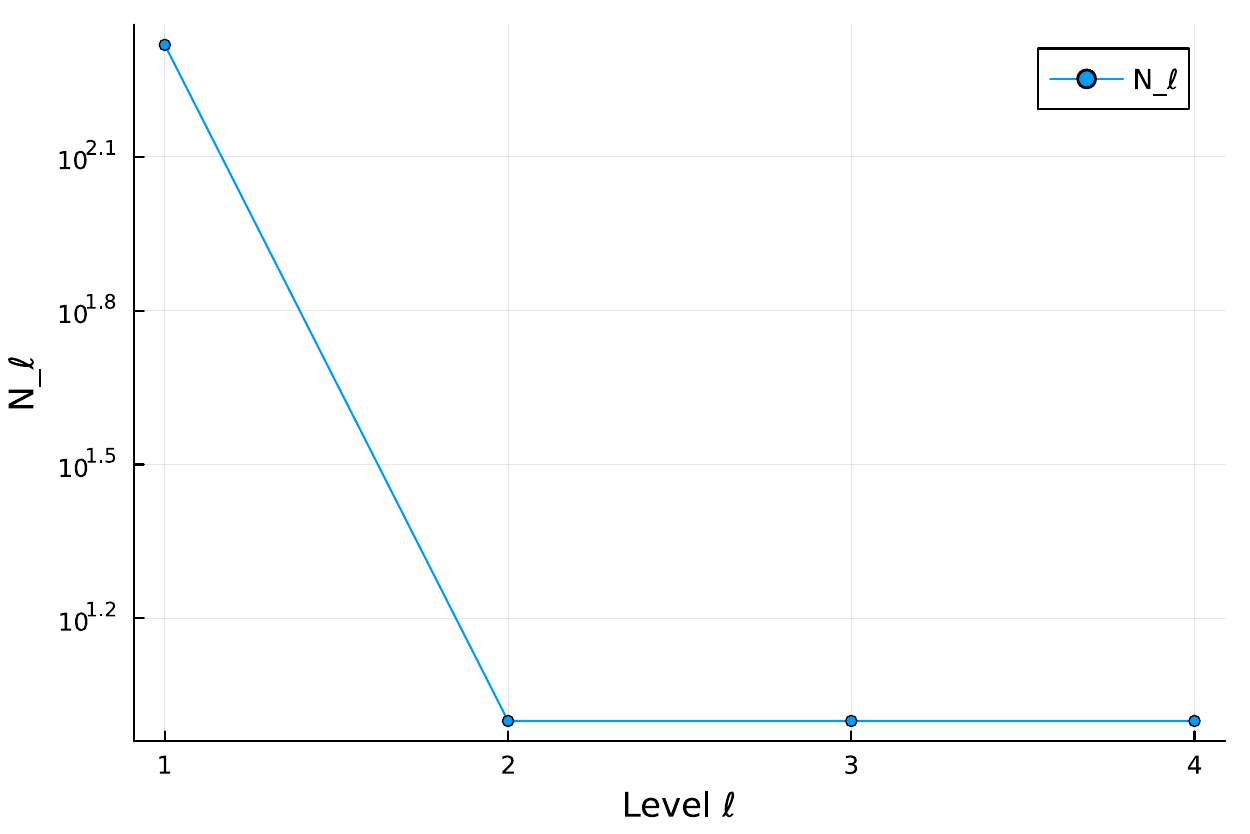}
			\caption*{MLMC Number of Realizations}
		\end{subfigure}
		\begin{subfigure}[b]{0.495\textwidth}
			\centering
			\includegraphics[width=\textwidth]{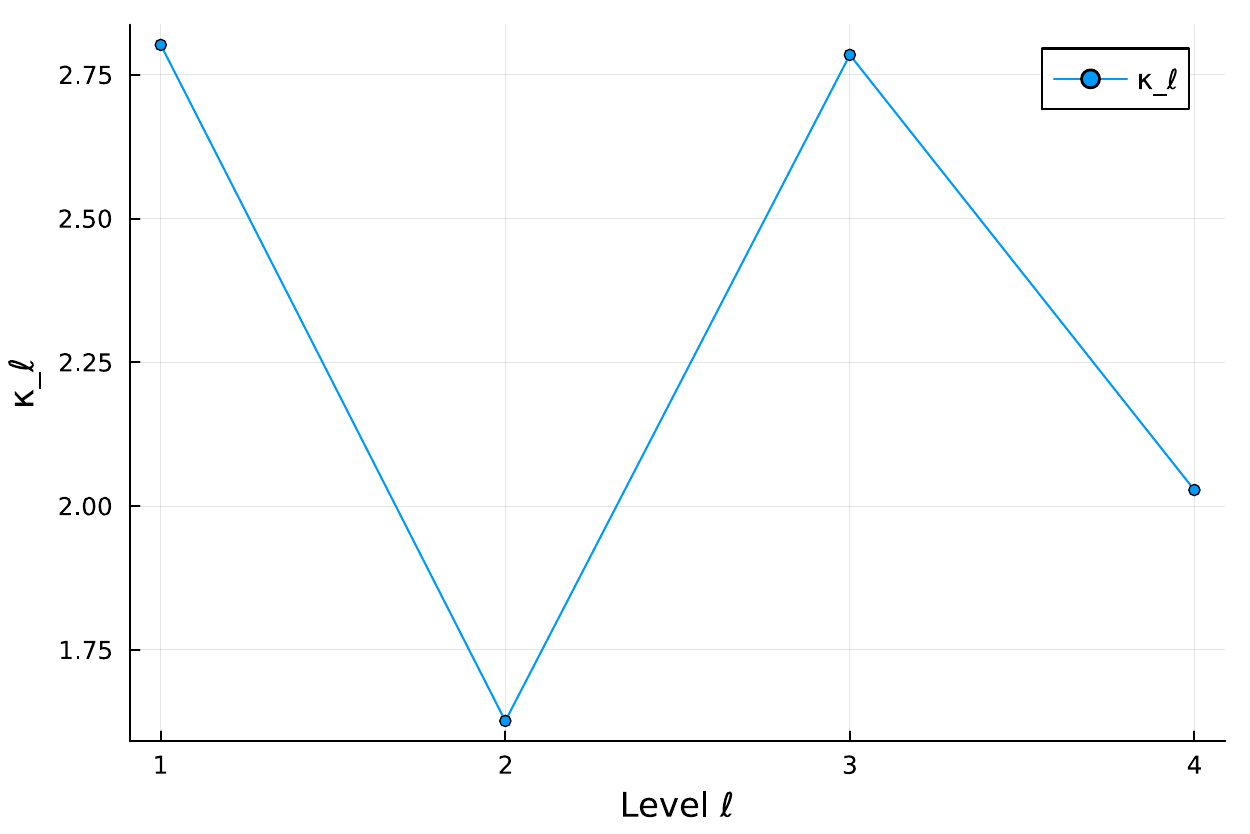}
			\caption*{MLMC Kurtosis}
		\end{subfigure}
		\begin{subfigure}[b]{0.495\textwidth}
			\centering
			\includegraphics[width=\textwidth]{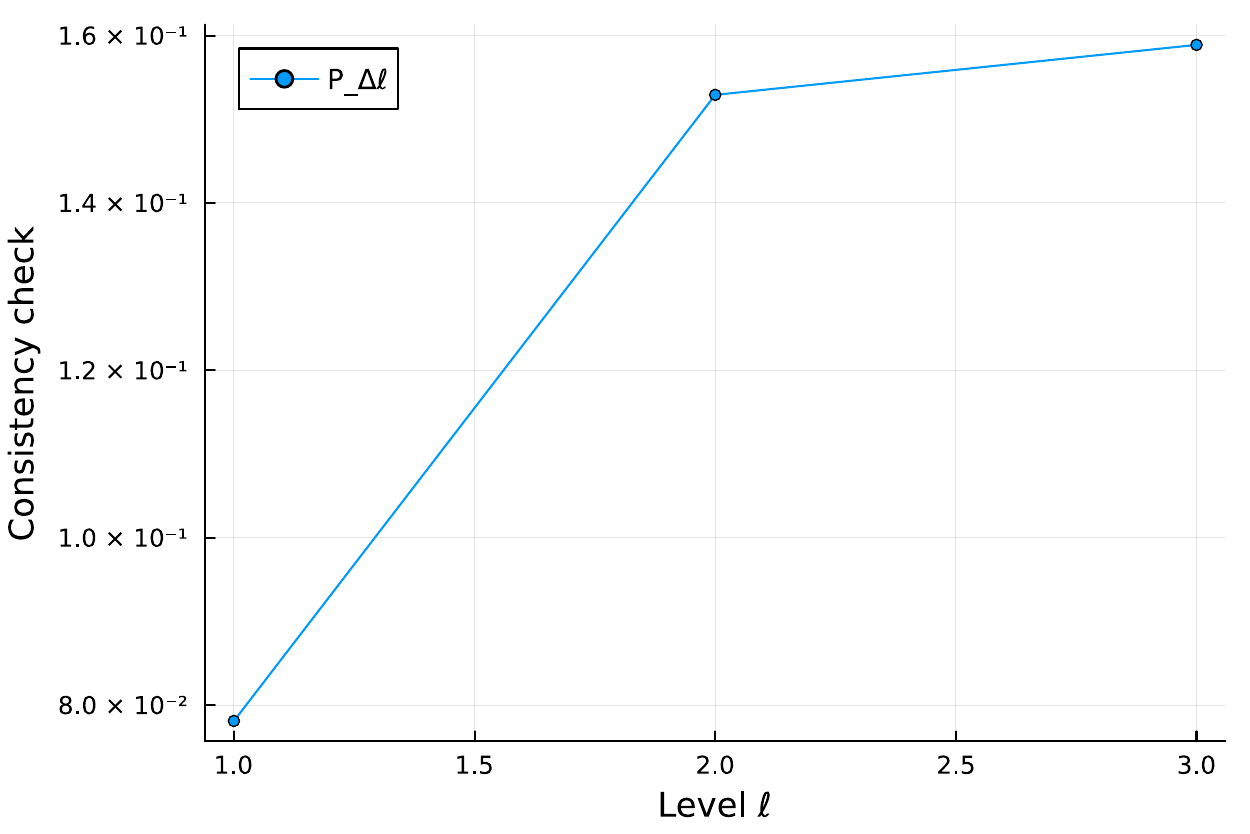}
			\caption*{MLMC consistency check}
		\end{subfigure}
		\caption{MLMC functional results for $c_2 = 0.9$, $K = 10^4$ particle histories, and $\epsilon = 1 \cdot 10^{-3}$}
		\label{fig:mlmc_results_0.9}
	\end{figure*}

	Figures \ref{fig:mlmc_results_0.5} presents the data on convergence of the MLMC algorithm in cases the functional defined over the whole domain for $c_2=0.5$, $K_{\ell}=10^4$, and $\epsilon = 1 \cdot 10^{-3}$.
	We see a monotonic decrease in variance and mean value of the correction factor with the addition of more levels. 
	In addition, we also see the cost estimate increases for a sample as we refine the computational grid.
	The results also show acceptable results for the consistency check and that the Kurtosis $\kappa_\ell$ is small relative to the number of initial samples.
	Definitions for the consistency check and Kurtosis are provided in Chapter 3 of \cite{giles-AN-2015}.
	Kurtosis demonstrates the convergence of the variance by giving an order of the number of samples needed for convergence, i.e. $N_\ell = \mathcal{O}(\kappa)$.
	The consistency check ensure that $P_{l-1,n}$ is being calculated in a way that ensures the validity of the telescoping summation.
	This checks consists of calculating:
	\vspace{-0.25cm}
	\[CC = \frac{ \big<P_{\ell-1}\big> - \big<P_{\ell}\big> + \big<\Delta P_{\ell}\big>}{3\Big(\sqrt{\mathbb{V}(P_{\ell-1})} + \sqrt{\mathbb{V}(P_{\ell})}+ \sqrt{\mathbb{V}(\Delta P_\ell)}\Big)}.
		\vspace{-0.25cm}
	\]
	$CC$ should be less than $1$ otherwise we have implemented the method incorrectly.
		In Figure \ref{fig:mlmc_results_0.9}, analysis for the whole domain functional with $c_2=0.9$, $K_{\ell}=10^4$, and $\epsilon = 1 \cdot 10^{-3}$ is shown.
	The behavior for $c_2 = 0.9$ is mostly the same as the previous case, except more samples were requested on the coarsest level.
\clearpage

	In Tables \ref{tab:mlmc_mixed_1e3} and \ref{tab:mlmc_mixed_1e4} we examine the behavior of the method using the whole domain functional for a wider range of $\epsilon$ values.
	The overall trend is as one would expect, as $\epsilon$ decreases the required number of samples requested increases.
	In addition, the samples that are requested are also placed on the coarsest grid as one would expect based upon Theorem 1.
	We find that the method passes the weak convergence check which means we have ran sufficient number of levels.
\begin{table}[htb]
	\centering
	\vspace*{-0.1cm}
	\caption{MLMC Linear Fit and Distribution of Samples for $K_{\ell} = 10^3$, 2 Materials}
	{\small
		\begin{tabular}{|c|c|c|c|c|c|c|c|c|c|}
			\hline
			$c_2$ & $\epsilon$       & $\alpha$ & $\beta$ & $\gamma$ & $N_0$ & $N_1$ & $N_2$ & $N_3$ & $\max{W_{\hat{\ell}}}$\\
			\hline
			0.1& $1 \cdot 10^{-2}$ & 2.01  & 2.06 & 0.72 & 10   & 10 & 10 & 10 & $3.3 \cdot 10^{-4}$ \\
			0.1& $5 \cdot 10^{-3}$ & 2.03  & 2.34 & 0.67 & 37   & 10 & 10 & 10 & $3.1 \cdot 10^{-4}$ \\
			0.1& $1 \cdot 10^{-3}$ & 2.04  & 2.30 & 0.66 & 663  & 10 & 10 & 10 & $3.3 \cdot 10^{-4}$ \\
			0.5& $1 \cdot 10^{-2}$ & 1.99  & 1.26 & 0.69 & 15   & 10 & 10 & 10 & $3.6 \cdot 10^{-4}$ \\
			0.5& $5 \cdot 10^{-3}$ & 2.09  & 1.97 & 0.68 & 32   & 10 & 10 & 10 & $3.3 \cdot 10^{-4}$ \\
			0.5& $1 \cdot 10^{-3}$ & 1.97  & 1.83 & 0.70 & 865  & 10 & 10 & 10 & $3.6 \cdot 10^{-4}$ \\
			0.9& $1 \cdot 10^{-2}$ & 1.97  & 1.60 & 0.72 & 24   & 10 & 10 & 10 & $4.9 \cdot 10^{-4}$ \\
			0.9& $5 \cdot 10^{-3}$ & 1.99  & 1.99 & 0.72 & 71   & 10 & 10 & 10 & $4.8 \cdot 10^{-4}$ \\
			0.9& $1 \cdot 10^{-3}$ & 1.96  & 2.85 & 0.62 & 1756 & 10 & 10 & 10 & $5.0 \cdot 10^{-4}$ \\
			\hline
		\end{tabular}
	}
	\label{tab:mlmc_mixed_1e3}
\end{table}

\begin{table}[htb]
	\centering
	\vspace*{-0.1cm}
	\caption{MLMC Linear Fit and Distribution of Samples for $K_{\ell} = 10^4$, 2 Materials}
	{\small
		\begin{tabular}{|c|c|c|c|c|c|c|c|c|c|}
			\hline
			$c_2$ & $\epsilon$       & $\alpha$ & $\beta$ & $\gamma$ & $N_0$ & $N_1$ & $N_2$ & $N_3$ & $\max{W_{\hat{\ell}}}$ \\
			\hline
			0.1& $1 \cdot 10^{-2}$ & 2.01  & 3.52 & 0.80 & 10   & 10 & 10 & 10 & $3.2 \cdot 10^{-4}$ \\
			0.1& $5 \cdot 10^{-3}$ & 1.99  & 3.50 & 0.80 & 10   & 10 & 10 & 10 & $3.3 \cdot 10^{-4}$ \\
			0.1& $1 \cdot 10^{-3}$ & 2.00  & 2.72 & 0.72 & 76   & 10 & 10 & 10 & $3.3 \cdot 10^{-4}$ \\
			0.5& $1 \cdot 10^{-2}$ & 2.00  & 2.56 & 0.78 & 10   & 10 & 10 & 10 & $3.6 \cdot 10^{-4}$ \\
			0.5& $5 \cdot 10^{-3}$ & 2.00  & 3.14 & 0.83 & 10   & 10 & 10 & 10 & $3.6 \cdot 10^{-4}$ \\
			0.5& $1 \cdot 10^{-3}$ & 2.00  & 2.62 & 0.82 & 63   & 10 & 10 & 10 & $3.6 \cdot 10^{-4}$ \\
			0.9& $1 \cdot 10^{-2}$ & 1.99  & 3.06 & 0.66 & 10   & 10 & 10 & 10 & $4.9 \cdot 10^{-4}$ \\
			0.9& $5 \cdot 10^{-3}$ & 2.00  & 3.48 & 0.73 & 10   & 10 & 10 & 10 & $4.9 \cdot 10^{-4}$ \\
			0.9& $1 \cdot 10^{-3}$ & 2.00  & 3.06 & 0.83 & 208  & 10 & 10 & 10 & $4.9 \cdot 10^{-4}$ \\
			\hline
		\end{tabular}
	}
	\label{tab:mlmc_mixed_1e4}
\end{table}

Table \ref{tab:mlmc_1e4_cell_8} shows results for $K_{\ell}=10^4$ histories and the coarse cell form of the functional for cell $8$.
$g_8^0$ is an interface cell between the two materials.
We see that the estimates of the linear fit factors is noisier than the functional based upon the integral over the whole domain.
In addition, fewer samples were requested for this definition of the functional.
This is likely due to the functional having a smaller value yielding a smaller estimate of the variance.
\begin{table}[htb]
	\centering
	\caption{MLMC Linear Fit Results and Distribution of Samples for $K_{\ell} = 10^4$ and $g_8^0$, 2 Materials}
	{\small
		\begin{tabular}{|c|c|c|c|c|c|c|c|c|c|}
			\hline
			$c_2$ & $\epsilon$       & $\alpha$ & $\beta$ & $\gamma$ & $N_0$ & $N_1$ & $N_2$ & $N_3$ & $\max{W_{\hat{\ell}}}$\\
			\hline
			0.1& $1 \cdot 10^{-2}$ & 2.00  & 3.12 & 0.78 & 10 & 10 & 10 & 10 & $3.5 \cdot 10^{-5}$ \\
			0.1& $5 \cdot 10^{-3}$ & 2.02  & 2.77 & 0.83 & 10 & 10 & 10 & 10 & $3.4 \cdot 10^{-6}$ \\
			0.1& $1 \cdot 10^{-3}$ & 2.02  & 2.68 & 0.92 & 10 & 10 & 10 & 10 & $3.4 \cdot 10^{-6}$ \\
			0.5& $1 \cdot 10^{-2}$ & 2.00  & 3.03 & 0.81 & 10 & 10 & 10 & 10 & $3.1 \cdot 10^{-5}$ \\
			0.5& $5 \cdot 10^{-3}$ & 2.03  & 2.33 & 0.78 & 10 & 10 & 10 & 10 & $3.1 \cdot 10^{-5}$ \\
			0.5& $1 \cdot 10^{-3}$ & 2.00  & 3.47 & 0.85 & 10 & 10 & 10 & 10 & $3.1 \cdot 10^{-5}$ \\
			0.9& $1 \cdot 10^{-2}$ & 2.01  & 1.63 & 0.79 & 10 & 10 & 10 & 10 & $3.0 \cdot 10^{-5}$ \\
			0.9& $5 \cdot 10^{-3}$ & 2.01  & 2.67 & 0.81 & 10 & 10 & 10 & 10 & $3.0 \cdot 10^{-5}$ \\
			0.9& $1 \cdot 10^{-3}$ & 1.98  & 2.26 & 0.82 & 10 & 10 & 10 & 10 & $3.1 \cdot 10^{-5}$ \\
			\hline
		\end{tabular}
	}
	\label{tab:mlmc_1e4_cell_8}
\end{table}

We also consider a 3 region test.
Region 1 is defined $0<x<0.25$, $c_1 = 0.9$, Region 2 is defined $0.25<x<0.75$, $c_2 = 0.5$, and Region 3 is defined $0.75<x<1.0$, $c_3 = 0.1$.
The total cross-sections we considered were $\Sigma_t = 1.0, 5.0$.
Results from this test are in Table \ref{tab:mlmc_mixed_3_materials_1e3}.
The results for $\Sigma_t = 1.0$ are similar to the previous test results while $\Sigma_t = 5.0$ requested fewer samples.
We analyzed results for $K_{\ell} = 10^4$ which are similar to the results we have presented in this paper: the requested number of samples is smaller than the $K_{\ell}=10^3$ case but the majority of the work is placed on the coarsest level.
\begin{table}[htb]
	\centering
	\caption{MLMC Linear Fit and Distribution of Samples for $K_{\ell} = 10^3$, 3 Materials}
	{\small
		\begin{tabular}{|c|c|c|c|c|c|c|c|c|c|}
			\hline
			$\Sigma_t$ &  $\epsilon$       & $\alpha$ & $\beta$ & $\gamma$ & $N_0$ & $N_1$ & $N_2$ & $N_3$ & $\max{W_{\hat{\ell}}}$\\
			\hline
			$1.0$      & $1 \cdot 10^{-2}$ & 1.98  & 2.73 & 0.70 & 10   & 10 & 10 & 10 & $2.5 \cdot 10^{-4}$ \\
			$1.0$      & $5 \cdot 10^{-3}$ & 2.02  & 3.28 & 0.67 & 30   & 10 & 10 & 10 & $2.6 \cdot 10^{-4}$ \\
			$1.0$      & $1 \cdot 10^{-3}$ & 2.07  & 1.65 & 0.65 & 491  & 10 & 10 & 10 & $2.5 \cdot 10^{-4}$ \\
			$5.0$      & $1 \cdot 10^{-2}$ & 1.98  & 2.50 & 0.61 & 10   & 10 & 10 & 10 & $5.3 \cdot 10^{-4}$ \\
			$5.0$      & $5 \cdot 10^{-3}$ & 2.03  & 2.99 & 0.71 & 10   & 10 & 10 & 10 & $5.0 \cdot 10^{-4}$ \\
			$5.0$      & $1 \cdot 10^{-3}$ & 1.98  & 2.36 & 0.69 & 25   & 10 & 10 & 10 & $5.2 \cdot 10^{-4}$ \\
			\hline
		\end{tabular}
	}
	\label{tab:mlmc_mixed_3_materials_1e3}
\end{table}

\section{Conclusions} \label{sec:con}

We have formulated a MLMC algorithm that uses a hybrid MC simulation to generate a scalar flux sample. 
In this study, we have analyzed the convergence of some functionals of the solution computed with the MLMC algorithm.
Results shown that we should expect $\mathcal{O}(\epsilon^2)$ computational complexity for this class of problems, as $\beta > \gamma$ for the cases we analyzed.
Our methodology has some unique features that warrant further discussion: in the original method the cost of generating a MLMC sample is driven primarily by the cost of solving the system of equations after a random sample is generated.
In our case, the main computational cost in the generation of a sample is the estimation of the Eddington and boundary factors using a MC simulation.
This is why our estimate of $\gamma$ is small when one would expect it to be closer to $1$ when we double the number of computational cells.
Estimates of $\beta$ are noisy, suggesting that outlier samples were affecting the results in some of the MLMC simulations.
Notable, our estimate of $\alpha$ is consistent across simulations at a value of $\alpha = 2$, as one would expect for a second-order discretization scheme like the one we used for the hybrid solve.

Overall, we see that the proposed MLMC algorithm converges the functional as we add computational levels.
Future work considerations is expanding the results to multi-dimensional Hybrid MC.
An important item of further research is formulating a method for estimation of correction relative to the neighboring level instead of direct calculation of $\big<P_{\ell}\big>$.
Formulation of another type of hybrid low-order equations can help reduce noise in the Hybrid MC samples that we use as input to the MLMC algorithm.

\section*{Acknowledgements}

The work of the first author (VNN) was supported under an University Nuclear Leadership Program Graduate Fellowship. Any opinions, findings, conclusions or recommendations expressed in this publication are those of the author(s) and do not necessarily reflect the views of the Department of Energy Office of Nuclear Energy. 
This work was also supported by the Center for Exascale Monte-Carlo Neutron Transport (CEMeNT), a PSAAP-III project funded by the Department of Energy, grant number DE-NA003967.

\bibliographystyle{elsarticle-num}
\bibliography{main}

\end{document}